\documentclass{amsart}

\usepackage{euscript}
\usepackage{epsfig}

\theoremstyle{plain}

\theoremstyle{definition}

\def\C{\mathbb C}

\def\Z{\mathbb Z}

\newcommand{\N}{\mathbb Z_{>0}}

\newcommand{\cE}{\mathcal E}
\newcommand{\cF}{\mathcal F}
\newcommand{\cG}{\mathcal G}

\newcommand{\cN}{\mathcal N}
\newcommand{\cO}{\mathcal O}

\newcommand{\cV}{\mathcal V}

\renewcommand{\t}{\widetilde}

\newcommand{\tX}{\widetilde X}

\DeclareMathOperator{\LF}{LF} 
\DeclareMathOperator{\di}{div} \DeclareMathOperator{\Ker}{Ker}
\DeclareMathOperator{\projan}{Projan}
\newcommand{\X}{(X,o)}

\newcommand{\defset}[2]{{\left\{#1\,\left| \,#2 \right. \right\}}}
\newcommand{\E}[1]{E_{v_{#1}}}
\newcommand{\thmref}[1]{Theorem~\ref{#1}}

\def\ten{\circle*{0.25}}
\renewcommand{\:}{\colon}

\begin{document}
\title{On the Casson Invariant Conjecture of Neumann--Wahl}
\author{Andr\'as N\'emethi}
\address{Department of Mathematics\\Ohio State University\\Columbus, OH 43210;
and}
\address{R\'enyi Institute of Mathematics\\Budapest, Hungary}
\email{nemethi@math.ohio-state.edu; \ nemethi@renyi.hu}
\author{Tomohiro Okuma}
\address{Department of Education, Yamagata University,
 Yamagata 990-8560, Japan.}
\email{okuma@e.yamagata-u.ac.jp}
\thanks{The first author is partially supported by NSF grant DMS-0605323, 
Marie Curie grant and OTKA
grants; the second author by the Grant-in-Aid for Young Scientists
(B), The Ministry of Education, Culture, Sports, Science and
Technology, Japan. }



\begin{abstract}
In the article we prove the Casson Invariant Conjecture of
Neumann--Wahl for splice type surface singularities. Namely, for
such an isolated complete intersection,  whose link is an integral
homology sphere, we show that the Casson invariant of the link is
one-eighth the signature of the Milnor fiber.
 \end{abstract}

\maketitle
\pagestyle{myheadings}

\markboth{{\normalsize A. N\'emethi and T. Okuma}}{ {\normalsize
On the Casson Invariant Conjecture}}

\newcommand{\bz}{{\mathbb Z}}
\newcommand{\zn}{Z_{min}}
\newcommand{\zmin}{Z_{min}}
\newcommand{\zk}{Z_K}
\newcommand{\co}{{\mathcal O}}
\newcommand{\x}{$(X,0)$\ }
\newcommand{\cl}{{\mathcal L}}
\newcommand{\no}{\noindent}
\newcommand{\bfc}{{\mathbb C}}
\newcommand{\bfq}{{\mathbb Q}}
\newcommand{\cale}{{\mathcal E}}
\newcommand{\calw}{{\mathcal W}}
\newcommand{\calv}{{\mathcal V}}
\newcommand{\cala}{{\mathcal A}}
\newcommand{\calr}{{\mathcal R}}
\newcommand{\calp}{{\mathcal P}}
\newcommand{\bc}{{\mathbb C}}
\newcommand{\br}{{\mathbb R}}
\newcommand{\bq}{{\mathbb Q}}
\newcommand{\vs}{\vspace{3mm}}
\newcommand{\si}{\sigma}

{\small

\section{Introduction}

\no Almost twenty years ago, Neumann and Wahl formulated the
following conjecture:\\

\no {\bf Casson Invariant Conjecture.} \cite{NW} \ {\em Let
$(X,o)$ be an isolated complete intersection surface singularity
whose link $\Sigma$ is an integral homology 3--sphere
Then the Casson invariant $\lambda(\Sigma)$ of the link is one-eighth the
signature of the Milnor fiber of $(X,o)$.}\\

The conjecture can be reformulated in terms of the geometric genus
$p_g$ of $(X,o)$ as well (see below in \ref{4.1}).

The conjecture is true for Brieskorn hypersurface singularities by a result 
of Fintushel and Stern \cite{FS}. This and additivity properties 
(with respect to splice decomposition) lead to the verification 
of the conjecture  for Brieskorn complete intersections, done independently
by Neumann--Wahl  \cite{NW} and Fukuhara--Matsumoto--Sakamoto \cite{FMS}.
For suspension hypersurface singularities it was verified in \cite{NW}. 
Some iterative generalizations,
related with cyclic coverings and using techniques of equivariant
Casson invariant and gauge theory, were covered by Collin and
Saveliev (cf. \cite{Co,CoS,CoS2}).

Recently, Neumann and Wahl have introduced an important family of
complete intersection surface singularities, the  {\em splice type
singularities}. In  \cite{NWuj} they treated the case when the
link is an integral homology sphere (the reader may consult in
\cite{NWnew,NWnew2,NWuj2} the case of rational homology sphere
links too). In \cite{NWuj}, they have also verified the above
conjecture (by a direct computation of the geometric genus) for
special splice type singularities (when the nodes of
the splice diagram `are in a line').

The goal of the present article is to verify the conjecture for an
arbitrary splice type singularity:\\

\no {\bf Theorem.} {\em The Casson Invariant Conjecture is true
for any splice type singularity with integral homology sphere link.}\\

This theorem can also be judged in the light of the following
expectation/conjecture (cf. \cite{NWuj,NWuj2}, see also section 6
in \cite{Stev}): any complete intersection surface singularity
with integral homology sphere link is of splice type.

The article is organized as follows. Section 2 is a review of
splice type surface singularities. The proof of the theorem is
broken into two parts: section 3 contains the analytic part, while
the topological/combinatorial part is in the last section 4.

The proof contains an inductive formula for the geometric genus.
The inductive step corresponds to the splice decomposition of the
link $\Sigma$ (cf. \cite{EN}). But we wish to emphasize that the
study of any {\em analytic} invariant with respect to the splice
decomposition is rather delicate since the splice decomposition
(gluing), rewritten in the language of plumbing, contains a purely
topological step (namely the 0--absorption, cf. \cite{EN}, \S
22) which cannot be represented in the world of negative definite
plumbing graphs. A new construction of section 4 surmounts this difficulty.

\vspace{2mm}

\no {\bf Acknowledgements.} \ The second  author thanks Marie
Curie Fellowship for Transfer of Knowledge for supporting his
visit at the R\'enyi Institute of Mathematics, Budapest, Hungary,
where the work on this paper was done. He is also grateful to
the members of the R\'enyi Institute for the warm hospitality.

\section{Splice type surface singularities (review)}

This section contains a brief introduction of {\em splice type
singularities}  in terms of ``monomial cycles'' (cf. \cite[\S
3]{Ouac-c}, \cite[\S 13]{NWuj2}). We always assume that all the
links are integral homology spheres.

The splice type singularities, which generalize Brieskorn complete
intersections,  were  introduced by Neumann and Wahl (see
\cite{NWnew}, \cite{NWuj}, \cite{NWuj2}). To any resolution (or,
negative definite plumbing) graph, associated with a fixed plumbed
3--manifold, one first associates a new weighted tree, called the
``splice diagram". From this one writes down the system of
(leading forms of) splice diagram equations. They define a splice
type singularity whose link is the original 3-manifold. (In the
next presentation the splice diagram is almost `hidden', the
interested reader may consult the above references for more
details.)

\subsection{Basic notations}\label{2.1}
Let $\X$ be a germ of a normal complex surface singularity and
$\pi\:\tX \to X$ a good resolution, i.e., its exceptional divisor
$E:=\pi^{-1}(o)$ has only simple normal crossings. Let $\{E_v\}_{v
\in \cV}$ denote the set of irreducible components of $E$. The
link $\Sigma$ of the singularity $X$ is an integral homology
sphere if and only if $E$ is a tree of rational curves and the
 intersection matrix $I(E):=(E_v\cdot E_w)_{v,w \in \cV}$ is unimodular.
Let $(m_{vw})_{v,w \in \cV}=-I(E)^{-1}$. Then every $m_{vw}$ is a
positive integer. We call an element of a group $\sum _{v\in
\cV}\Z E_v$  a cycle. For any cycle $D=\sum _{w \in \cV}a_wE_w$,
we write $m_w(D)=a_w$.

Let $\delta_v=(E-E_v)\cdot E_v$, the number of irreducible
components of $E$ intersecting $E_v$. A curve $E_v$ (or its index
$v$) is called an {\itshape end} (resp. a {\itshape node}) if
$\delta_v=1$ (resp. $\delta_v\ge 3$). We denote by $\cE$ (resp.
$\cN$) the set of indices of ends (resp. nodes). A connected
component of $E-E_v$ is called a {\itshape branch} of $v$.

\subsection{Monomial cycle}
For any $v \in \cV$, let $E^*_v=\sum_{w \in \cV}m_{vw}E_w$. Then
$E^*_v\cdot E_w=-\delta_{vw}$ for every $w \in \cV$,  where
$\delta _{vw}$ denotes the Kronecker delta.
An element of a semigroup $\sum _{w \in \cE}\Z_{\ge 0}E^*_w$,
 where $\Z_{\ge 0}$ is the set of nonnegative integers, is
 called a {\itshape monomial cycle}.
Let $\C\{z\}:=\C\{z_w ; w \in \cE\}$ be the convergent power
series ring in $\# \cE$  variables. Then  for a monomial cycle
$D=\sum _{w \in \cE}\alpha_wE_w^*$, we
 associate a monomial  $z(D):=\prod_{w\in \cE} z_w^{\alpha_w} \in
 \C\{z\}$.

\subsection{Degree and order associated with $v$.} We fix $v \in \cV$.
One defines the $v$-weight of any variable $z_w$, $w \in \cE$, by
$m_{vw}$. Note that $\gcd\{m_{vw}\}_{w \in \cV}=1$. If $D$ is a
monomial cycle, then the  $v$-degree of $z(D)$ is equal to
$m_v(D)$.

For any  $f\in \C\{z\}$ write $f=f_0+f_1 \in \C\{z\}$, where $f_0$
is a nonzero quasihomogeneous polynomial with respect to the
$v$-weight and  $f_1$ is a series in monomials of higher
$v$-degrees. Then we call $f_0$ the {\itshape $v$-leading form} of
$f$, and denote it by $\LF_{v}(f)$.
We define the $v$-order of $f$ by $v\mbox{-ord}(f):=v\mbox{-deg}(\LF_v(f))$.

\subsection{Monomial Condition}\label{ss:mc}
We say that $E$ satisfies the monomial condition if for  any node
$v$ and any  branch $C$ of $v$, there exists a monomial cycle $D$
such that $D-E^*_v$ is an effective cycle supported on $C$. In
this case, $z(D)$ is called an {\itshape admissible monomial}
belonging to the branch $C$. Note that a chain $C$  always admits
a monomial cycle.

The monomial condition is equivalent to the semigroup and
congruence conditions of Neumann-Wahl (see \cite[\S 13]{NWuj2})
(although the congruence condition is trivial in our case.)

\subsection{Neumann-Wahl system}
Assume that the monomial condition is satisfied. For any fixed
node $v$ let $C_1, \dots , C_{\delta_v}$ be the branches of $v$.
Suppose $\{m_1, \ldots ,m_{\delta_v}\}$  is a set of admissible
monomials such that $m_i$ belongs to $C_i$ for  $i=1, \ldots
,\delta_v$. Let $F=(c_{ij})$, $c_{ij} \in \C$, be any
$((\delta_v-2) \times \delta_v)$-matrix such that all the maximal
minors have rank $\delta_v-2$. We define polynomials $f_1,
\ldots ,f_{\delta_v-2}$ by
$$\begin{pmatrix}
 f_1\\ \vdots \\ f_{\delta_v-2}
\end{pmatrix}
= F
\begin{pmatrix}
 m_1\\ \vdots \\ m_{\delta_v}
\end{pmatrix}.
$$
We call the set $\{f_1 , \dots ,f_{\delta_v-2}\}$ a {\itshape
Neumann-Wahl system} at $v$. If we have a Neumann-Wahl system
$\cF_v$ at every node $v$ then we call the set $\cF:= \bigcup_{v
\in \cN}\cF_v$ a {\itshape Neumann-Wahl system} associated with
$E$. Note that $\# \cF=\#\cE-2$.

\subsection{Splice diagram equations}\label{ss:sde}
Consider a finite set of germs
$$\defset{f_{vj_v}}{v \in \cN, \; j_v=1, \dots
,\delta_v-2}\subset \C\{z\}.$$ If the set
$\defset{\LF_v(f_{vj_v})}{v \in \cN, \; j_v=1, \dots ,\delta_v-2}$
is a Neumann-Wahl system associated with $E$  then a system of
equations
$$f_{vj_v}=0, \quad v \in \cN, \quad j_v=1, \dots ,\delta_v-2,$$
is called the {\itshape splice diagram equations} and the germs
$\{f_{vj_v}\}$ {\em splice diagram functions}. The germ of
singularity defined by the splice diagram equations in
$(\C^{\#\cE},o)$ is called a {\itshape splice type singularity}
(associated with the combinatorics of $E$).

We note that a  splice type  singularity is  an equisingular
deformation of a singularity defined by a  Neumann-Wahl system
associated with $E$ (cf. \cite[(4.)]{Ouac-c}).

In fact, all the above construction of the splice equations
(including the validity of the monomial condition) depends only on
the dual graph associated with $E$, hence only on a fixed plumbing
graph of the link $\Sigma$.

\subsection{Theorem}\label{t:2.1} \cite[(2.1)]{NWuj} {\em  Suppose that $(Z,o)$ is a splice type singularity
associated with $E$ (or, with a graph of $\Sigma$). Then  $(Z,o)$
is an isolated complete intersection surface singularity whose
link is $\Sigma$.}

\subsection{End-Curve Condition} We keep the notations of (\ref{2.1}).
We say that $\tX$ satisfies
the {\em end-curve condition} if for each $w \in \cE$ there exists an
irreducible curve $H_w \subset \tX$, not an exceptional curve, and
a function $f \in H^0(\cO_{\tX}(-E_w^*))$ such that
$\di(f)=E_w^*+H_w$. We call such $f$ an {\itshape end-curve
function} of $E_w$. The end-curve condition is equivalent to the
fact that for every $w \in \cE$ the linear system of
$\cO_{\tX}(-E_w^*)$ has no fixed component in $E$.

Notice that if $\tX$ satisfies the end-curve condition then so do
the minimal good resolution and any resolution obtained from $\tX$
by  blowing-up the singular points of $E$.

If $(X,o)$ is defined by splice diagram equations, then its
resolution $\tX$ satisfies the end-curve condition (since the
coordinate functions serve as end-curve functions). The converse
is guaranteed by the following:

\subsection{Theorem}\label{t:nw-end-curve} \cite[(4.1)]{NWuj} \ 
{\em Let $\pi:\tX\to X$ be a good resolution of a singularity $(X,o)$.
 If $\tX$ satisfies the end-curve condition, then $E$
satisfies the monomial condition and  $X$ is of splice type. In
fact, if
$$ \psi\: \C\{z_w; w \in \cE\} \to \cO_{X,o}$$
is a homomorphism of $\C$-algebras which maps each $z_w=z(E_w^*)$
to an end-curve function of $E_w$, then $\psi$ is surjective and
$\Ker \psi$ is generated by splice diagram functions.}

\section{Filtration and $p_g$-formula}

Assume that $(X,o)$ is of splice type, hence its resolution
satisfies the end-curve condition. We also assume the presence of
a node in $E$ (otherwise the graph is $A_n$). The main result of
this section presents an inductive  formula  for the geometric
genus $p_g(X)$ of $X$ with respect to `cutting' $E$ in two parts
(corresponding to a certain `splice decomposition' of $\Sigma$).
It is convenient to do this cutting near a vertex $v$ 
which isolates either a `star-shaped subgraph' or a `string'. 
The inductive formula involves a filtration of $\cO_{X,o}$
associated with $v$. It  can be proved by repeating the
line of arguments of \cite{Opg} for our slightly more general
situation (in \cite{Opg} $v$ is node); for the convenience of the
reader we provide the main points of the proof emphasizing the
differences.

Let $v_1$ be an end-node, i.e., $\E1$  is an end of the minimal
reduced connected cycle containing all nodes of $E$. Let $v_2$ be
the node which is nearest to $v_1$. We fix  $v\in\cV$ in that
branch of $v_2$ which contains $v_1$ such that $v$ is not an end.
(If $v_1$ is the only node, let $v$ be any vertex which is not an
end.) Let $C_1, \dots ,C_{\delta_v}$ denote the branches of $E_v$.

\subsection{Filtration with respect to $E_v$}
The $v$-order defines a filtration $\{F_n\}$ of $\C\{z\}$ by
$$F_n=\{f\in\C\{z\}\ |\  v\mbox{-ord}(f)\geq n\} \ \ (n\in \Z_{\ge 0}).$$
This induces a filtration of ideals  $I_n\subset \cO_{X,o}$ 
($n \in \Z_{\ge 0}$) by $I_n:=\psi(F_n)$, where 
$\psi:\C\{z\}\to \cO_{X,o}$ is the  natural projection (\ref{t:nw-end-curve}).
Let  $\cG$ denote the associated graded algebra
$\bigoplus_{n\ge0}\cG_n$, where $\cG_n=I_n/I_{n+1}$.

\subsection{Proposition}\label{p:graded} (Cf. \cite[(2.6)]{NWuj2}.) \ {\em  Let $\{f_{wj_w}\}\subset \C\{z\}$ be a set of
splice diagram functions as in (\ref{ss:sde}), and $I$ the ideal
of the polynomial ring $\C[z]$ generated by $v$-leading forms
$\{\LF_v(f_{wj_w})\}$. Then $\cG\cong \C[z]/ I$ and it is a
reduced complete intersection ring.}

\begin{proof}
The statement can be reduced to the case when $v$ is a node
considered in \cite[(2.6)]{NWuj2} (or \cite[(4.4)]{Ouac-c}) as
follows. Let  $D_1,\ldots, D_{\delta_{v_1}}$ denote the branches
of $v_1$; suppose that $D_j$ is a chain of curves and the variable
$z_j$ corresponds to the end of $D_j$ for  $1\le j \le
\delta_{v_1}-1$. If $E_v$ is on one of these branches, let $D_1$ be that 
branch.    Then we may assume that the Neumann-Wahl system
satisfies the following properties:
 \begin{align*}
\LF_w (f_{wj_w}) &=
m_{wj_w}+a_{wj_w}m_{w\delta_w-1}+b_{wj_w}m_{w\delta_w},  &
&w\in\cN, \; 1\le j_w \le \delta_w-2, \\
 m_{v_1j} & =  z_j^{\alpha_j}, & & 1\le j \le \delta_{v_1}-1.
\end{align*}
Note that  $\LF_{v}  (f_{wj_w})=\LF_v(\LF_{w} (f_{wj_w}))$ by
\cite[(3.8)]{Ouac-c}; while one gets $\LF_v$ of $\LF_{w}
(f_{wj_w})$ by deleting the monomials belonging to a branch
containing $v$. Suppose the admissible monomial
$M_{v_1\delta_{v_1}}$ contains $z_t$. Then we get that
$$ \{\LF_v(f_{wj_w})\}\bigcup \{z_1, z_t\} =
 \{\LF_{v_1}(f_{wj_w})\}\bigcup \{z_1, z_t\}.$$
In the above-mentioned articles it is shown that the second set is
a regular sequence and $\{\LF_{v_1}(f_{wj_w})\}\bigcup \{z_1\}$
defines a curve in $\C^{\# \cE}$  smooth off the origin.
Since $\{\LF_{v_1}(f_{wj_w})\}$ is also a regular sequence, 
 $\cG\cong \C[z]/ I$ by \cite[(3.3)]{NWuj2}. For the reducedness, see the
 argument at the end of section 3 of  \cite{NWuj2}.
\end{proof}

\subsection{Hilbert series}
Let $H(t)$ denote the Hilbert series of the graded ring $\cG$,
i.e.,
$$H(t)=\sum_{i\ge 0}\,(\dim \cG_i)\,t^i.$$
From (\ref{p:graded}), by a well--known formula, valid for complete 
intersections, one has:

\subsection{Proposition}\label{p:hp}
$$ H(t)=\prod_{w\in \cV} \left(1-t^{m_{vw}}\right)^{\delta_w-2}.$$

\no From (\ref{p:graded}) one may also read the $a$--invariant of
Goto and Watanabe associated with $\cG$ (namely
$\sum_{w\in\cN,j_w}v\mbox{-deg}(LF_v(f_{wj_w}))-
\sum_{w\in\cE}v\mbox{-deg}(z_w)$,
cf.  \cite[(3.1.4)]{GW}):
$$a(\cG)= \sum_{w\in\cV}(\delta_w-2)m_{vw}.$$

\subsection{}
Let $\pi'\: \tX \to X'$  be the morphism which contracts each
connected component of the  divisor $E-E_v$ to a normal point and $\pi_1\:
X'\to X$ the natural morphism. Let $E'=\pi'(E_v)$.

\subsection{Lemma}\label{l:filt1} {\em
\begin{enumerate}
\item There exists  $k \in \N$ such that $\cO_{\tX}(-kE_v^*)$ is $\pi$-generated.
\item  For any $n \in \Z_{\ge 0}$, $I_n=
H^0(\cO_{X'}(-nE'))
=H^0(\cO_{\tX}([-n\pi'^*E']))$. 
\item   $-E'$ is $\pi_1$-ample and $\pi_1$ coincides with the filtered blowing up
$$\projan_X\left(\bigoplus_{n\ge 0}\pi_{1*}\cO_{X'}(-nE')\right) \to X.$$
\item $H^1(\cO_{X'}(-nE'))=0$ for $n > a(\cG)$.
\end{enumerate}}

\begin{proof} First we treat (1). Let $C$ be one of the branches
of $v$, and $E_u$ be that component of $C$ which intersects $E_v$.
Fix any end  $E_w$ of $E$ in  $C$, and set $A:=
-\sum_{E_v\subset C}\,\det(-C)\,I(C)^{-1}_{vw}\,E_v$. This is an
effective cycle supported by $C$ such that $m_u(A)E_v^*+A$ is a
monomial cycle, namely $\det(-C)E_w^*$. Therefore, if we fix two
branches $C_i$ ($i=1,2$), we can construct two positive integers
$k_i$ and effective cycles $A_i$ supported by $C_i$ such that
$k_iE_v^*+A_i$ is monomial. Set $k=k_1k_2$ and $f_i \in
H^0(\cO_X)$ the image of the monomial $z(kE_v^*+k_jA_i)$, where
$i,j=1,2$ and $i\neq j$. Then it is easy to see that $f_1, f_2 \in
H^0\left(\cO_{\tX}(-kE_v^*)\right)$ generate $\cO_{\tX}(-kE_v^*)$.

The proof of (2)--(4) is same as the proof of 
\cite[(3.3--3.4)]{Opg} based on part (1) and  filtered ring theory \cite{To2}.
\end{proof}

Since $m_{vv}\pi'^*E'=E_v^*$ is integral, from (2) we get for any $m\in\Z_{>0}$
\begin{equation}\label{star}
I_{mm_{vv}}=H^0(\cO_{\tX}(-mE^*_v)). 
\end{equation}

\subsection{$p_g$-formula}
Let $(X_i, x_i)$ denote the singularity germ obtained by
contracting the branch $C_i$ (hence $X_i\subset X'$).  Fix 
$k\in \N$ such that $\pi'_*(kE_v^*)$ is Cartier. 
Its existence  is guaranteed by the proof of  
(\ref{l:filt1})(1) (namely by the existence of functions $f_i$). 

\subsection{Theorem}\label{t:pg} {\em Set
$P(n)=\sum_{i=0}^{n-1}(\dim \cG_i)$ for $n\in\N$. Then for any
integer $m > a(\cG)/km_{vv}$}
$$p_g(X,o)=P(mkm_{vv})-\chi(mkE_v^*) +\sum_{i=1}^{\delta_v}p_g(X_i,x_i).$$

\begin{proof}
We write $D=mkE_v^*$ and $D'=mkm_{vv}E'$. Consider the exact
sequence
$$0 \to \cO_{\tX}(-D) \to \cO_{\tX} \to \cO_{D}\to 0.$$
By equation (\ref{star}) one has
$$\dim_{\C}H^0(\cO_{\tX})/H^0(\cO_{\tX}(-D))=\dim_{\C}I_0/I_{mkm_{vv}}=
P(mkm_{vv}).$$
Therefore it suffices to show that
\begin{equation}\label{eq:nb}
h^1(\cO_{\tX}(-D))=\sum_{i=1}^{\delta_v}p_g(X_i,x_i).
\end{equation}
From the spectral sequence
$$E_2^{i,j} = H^i(R^j\pi'_*\cO_{\tX}(-D)) \Rightarrow
H^n(\cO_{\tX}(-D))$$ one gets the exact sequence
$$H^1(\pi'_*\cO_{\tX}(-D)) \longrightarrow
H^1(\cO_{\tX}(-D))\stackrel{\alpha}{\longrightarrow}
H^0(R^1\pi'_*\cO_{\tX}(-D)) \to 0.$$ 
Since $\pi'_*D=D'$ is Cartier and 
$\pi'^*E'=E_v^*/m_{vv}$, by projection formula
and (\ref{l:filt1})(4),
$$H^1(\pi'_*\cO_{\tX}(-D))=H^1(\pi'_*\cO_{\tX}(-\pi'^*D'))
=H^1(\cO_{X'}(-D'))=0,$$
hence $\alpha$ is an isomorphism. Clearly the support of
$R^1\pi'_*\cO_{\tX}(-D)$ is in the set $\{x_i\}_i$. Since
$D=\pi'^*D'$ and $D'$ is Cartier, we have
$$\left(R^1\pi'_*\cO_{\tX}(-D)\right)_{x_i}\cong
\left(R^1\pi'_*\cO_{\tX} \otimes\cO_{X'}(-D')\right)_{x_i} \cong
\left(R^1\pi'_*\cO_{\tX}\right)_{x_i}.$$ This proves
\eqref{eq:nb}.
\end{proof}

\subsection{Periodic constants}\label{pc}
In the formula of (\ref{t:pg}), $P(mkm_{vv})-\chi(mkE_v^*)$ is
independent of $m\gg 0$. Let us formulate this property. Let
$F(t)=\sum _{i\ge 0}a_it^i$ be a formal power series. Set
$P_F(n)=\sum_{i=0}^{n-1}a_i$. Suppose that $P_F(kn)$ is a
polynomial function of $n$ for some $k \in \N$. Then for any $k'
\in \N$ satisfying this property, the constant terms of $P_F(kn)$
and $P_F(k'n)$ are the same. We call this constant the {\itshape
periodic constant} of $F(t)$ and denote it by $F|_{pc}$. For
example, the formula of \thmref{t:pg} is
$$p_g(X,o)=H |_{pc}+\sum  _{i=1}^{\delta_v}p_g(X_i,x_i).$$
Note that if two formal power series $F_1(t)$ and $F_2(t)$ have
periodic constants then
$$(F_1+F_2)|_{pc}=F_1|_{pc}+F_2|_{pc}.$$

\subsection{Remarks}\label{r:Opg} (1) The message of the above theorem
(\ref{t:pg}) is the following. Fix a negative definite resolution graph
$\Gamma$ (or its plumbed 3-manifold $\Sigma$) 
which satisfies the monomial condition.  Then 
there exists a splice type singularity whose link is $\Sigma$ (via
Neumann-Wahl system/construction), and the geometric genus of
any such analytic structure is independent of the choice of the 
Neumann-Wahl system (by 4.3 of \cite{Ouac-c} and 10.1 of \cite{NWuj2}),
hence it depends only on the combinatorics of $\Gamma$ --- 
it will be denoted by $p_g(\Gamma)$. Moreover, the geometric genus
of any such analytic structure satisfies the inductive formula (\ref{t:pg}).

(2) The end-curve condition and splice diagram equations are
considered for any surface singularity with rational homology
sphere links. For this case, and under the assumption that $v$ is a node,
a $p_g$-formula is established in \cite{Opg}.

\section{The proof of the main theorem}

\subsection{}\label{4.1} In this section we prove the main
Theorem stated in the introduction. For this, we will apply
theorem (\ref{t:pg}) in several different situations.

We start with $(X,o)$, a splice type singularity (whose minimal
good resolution $\tX$ satisfies the end-curve condition), and
whose link $\Sigma$ is an integral homology sphere. $\Gamma$ will
denote the dual resolution graph associated with $\tX$; $\cV$
stays for the set of vertices of $\Gamma$.

It is well-known that in the conjecture one may replace the
signature $\sigma$ of the Milnor fiber by the geometric genus
$p_g$  of $(X,o)$ via a formula of Durfee and Laufer \cite{Du}
$\sigma+8p_g+c(\Gamma)=0$, where $c(\Gamma):=
K_{\tX}^2+b_2(\tX)$ (here $K$ denotes the canonical
class and $b_2$ the second Betti number) depends only on the
combinatorics of $\Gamma$. 

If $\Gamma$ has only one node then Theorem follows from 
\cite[(7.7)]{NWuj}. A different argument runs as follows: The above $p_g$
formula (namely $p_g=H|_{pc}$, cf. \ref{pc}) is valid for any splice type
analytic structure (supported by the same topological type), hence
we may assume that $(X,o)$ admits a good $\C^*$--action. Then
$(X,o)$ is automatically a Brieskorn complete intersection (cf.
\cite{Neu}), hence we may apply \cite{NW} too. Therefore, in the
sequel we assume that $\Gamma$  has at least two nodes.

\subsection{}\label{moncon}
In our inductive procedure the following fact is crucial
(cf. also with \ref{r:Opg}): for any
splice type singularity with (not necessarily minimal) fixed good
resolution graph $\Gamma$, consider any connected subgraph
$\Gamma' \subset \Gamma$. Let $E'$ denote the
reduced connected cycle corresponding to $\Gamma'$. Then a
neighborhood of $E'$ satisfies the end-curve condition (see
\cite[(2.15)]{Opg}); hence  $E'$ (or $\Gamma'$)
satisfies the monomial
condition by (\ref{t:nw-end-curve}). Therefore if $(X',x')$
denotes the normal surface singularity obtained by contracting
$E'$, then $p_g(X',x')$ is topological computable from $\Gamma'$.
We write $p_g(X',x')=p_g(\Gamma')$.

\subsection{Starting the induction}
Let $v_1 \in \cN$ be an end-node of $\Gamma$ and $v_2 \in \cN$ the
node which is nearest to $v_1$. Write  $\Gamma' \subset \Gamma$
for the  branch of $v_1$ containing $v_2$, and denote by $w$ that
vertex in $\Gamma'$  which is connected by $v_1$ in $\Gamma$. Let
$H_{\Gamma,v_1}(t)$ be the Hilbert series of the associated graded
ring of the filtration with respect to $E_{v_1}$. Since
$p_g(\Gamma'')=0$ for all the branches $\Gamma''$ of $v_1$ which are chains,
by (\ref{t:pg}) we have the following.

\subsection{Proposition}\label{p:1}
 $p_g(\Gamma)=H_{\Gamma,v_1}|_{pc}+p_g(\Gamma')$.

\subsection{}\label{FIG} One needs to modify this inductive step, since,
e.g., in general, $\Gamma'$ is not unimodular.
The next goal is to fit $\Gamma'$ into a `good' inductive
procedure which is compatible with the splice decomposition of
$\Sigma$ (see also (\ref{additivity}) for more motivation). 

Let $\Delta $ be the splice diagram of $\Sigma$. For the
correspondence between  $\Gamma$ and $\Delta$ see  \cite[\S 22]{EN}. We
decompose $\Delta$ as the splice of two splice diagrams $\Delta_1$
and $\Delta_2$, where $\Delta_i\ni v_i$.

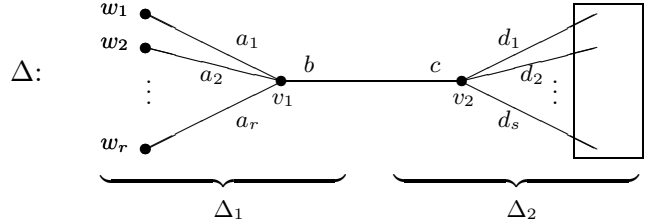
\begin{center}
\setlength{\unitlength}{0.6cm}
 \begin{picture}(12,5)(-6,-3)
\put(-2,0){\line(1,0){4}}
\put(-5,1.5){\ten}
\put(-5,0.75){\ten}
\put(-5,-1.5){\ten}
\put(-6,1.5){$w_1$}
\put(-6,0.75){$w_2$}
\put(-6,-1.5){$w_r$}
\put(-2,0){\ten}
\put(-2.2,-0.5){$v_1$}
\put(-1.5,0.2){$b$}
\put(2,0){\ten}
\put(1.8,-0.5){$v_2$}
\put(1.3,0.2){$c$}
\put(2,0){\line(2,1){3}}
\put(2,0){\line(4,1){3}}
\put(2,0){\line(2,-1){3}}
\put(-2,0){\line(-2,1){3}}
\put(-2,0){\line(-4,1){3}}
\put(-2,0){\line(-2,-1){3}}
\put(-5,1.5){\ten}
\put(-5,0.75){\ten}
\put(-5,-1.5){\ten}
\put(-6,1.5){$w_1$}
\put(-6,0.75){$w_2$}
\put(-6,-1.5){$w_r$}
\put(-5,-0.5){$\vdots$}
\put(-3,0.8){$a_1$}
\put(-3.8,0){$a_2$}
\put(-3,-1){$a_r$}
\put(4,-0.5){$\vdots$}
\put(2.8,0.8){$d_1$}
\put(3.3,0){$d_2$}
\put(2.8,-1){$d_s$}
\put(-6,-2){$\underbrace{\qquad\qquad\qquad\qquad\qquad}$}
\put(-3.5,-3){$\Delta_1$}
\put(0.5,-2){$\underbrace{\qquad\qquad\qquad\qquad\qquad}$}
\put(3,-3){$\Delta_2$}
\put(-8,0){{\large $\Delta$:}}
\put(4.5,-1.7){\framebox(1.5,3.4)}
 \end{picture}

\end{center}

\no
Let $\Gamma_i$ be the minimal plumbing graph associated with
$\Delta_i$ ($i=1,2$) (cf. \cite{EN}). Next we review how one can
recover from $\Delta_1$ and $\Delta_2$ the maximal chain
$\Gamma_{0}$ in $\Gamma$ which connects $v_1$ and $v_2$, and we
relate $\Gamma_2$ with $\Gamma'$. Set $a:=\prod_{i=1}^ra_i$ and
$d:=\prod_{i=1}^sd_i$. Express $a/b$ as a continued fraction
$[\alpha_0, \dots,\alpha_m]$,  i.e. 
$a/b=\alpha_0-1/(\alpha_1-1/(\ldots-1/\alpha_m))$ with 
$\alpha_i \in \Z_{>0}$ and $\alpha_i\ge 2$ for $i\ge 1$.
Similarly, write $d/c=[\beta_0, \dots ,\beta_n]$,
and consider the chain:

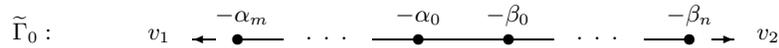
\begin{center}
\setlength{\unitlength}{0.6cm}

 \begin{picture}(12,1.5)(-2,-0.5)
\put(0,0){\line(1,0){1}} \put(1.5,0){. . .}
\put(3,0){\line(1,0){4}} \put(7.5,0){. . .}
\put(9,0){\line(1,0){1}} \multiput(0,0)(4,0){2}{\ten}
\multiput(6,0)(4,0){2}{\ten} \put(-0.5,0.4){$-\alpha_m$}
\put(3.5,0.4){$-\alpha_0$} \put(5.5,0.4){$-\beta_0$}
\put(9.5,0.4){$-\beta_n$} \put(-0.5,0){\vector(-1,0){0.5}}
\put(-2,0){$v_1$} \put(10.5,0){\vector(1,0){0.5}}
\put(11.5,0){$v_2$} \put(-5,0){$\t\Gamma_0:$}
 \end{picture}
\end{center}

\no This is not a `minimal' graph, and it is equivalent via
(topological) plumbing calculus with the (minimal) chain
$\Gamma_0$. This calculus runs as follows:
 by the positivity of the edge determinant
$bc>ad$ (cf. \cite[\S 1]{NWuj}), $\alpha_0=1$ or $\beta_0=1$. Then
one successively blows down the vertices whose weight are $-1$,
and at some moment inevitably a vertex with weight $0$ appears;
then  one makes a 0-absorption by the rule:
\setlength{\unitlength}{0.6cm}
 \begin{picture}(8,0.9)(0.5,0)
\put(0.3,0){\line(1,0){4}} \multiput(1,0)(1.5,0){3}{\ten}
\put(0.5,0.4){$-e_1$} \put(2.4,0.4){$0$} \put(3.5,0.4){$-e_2$}
\put(4.7,0){$\longrightarrow$} \put(6,0){\line(1,0){2.4}}
\put(7.2,0){\ten} \put(6.5,0.4){$-e_1-e_2$}
 \end{picture}\ .

\no In this way one  gets $\Gamma_0$. On the other hand, the
corresponding maximal chain  (with determinant $c$) in $\Gamma_2$ has the
following form:

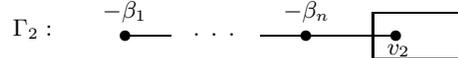
\begin{center}
\setlength{\unitlength}{0.6cm}
 \begin{picture}(7,1.5)(0,-0.5)
\put(0,0){\line(1,0){1}} \put(1.5,0){. . .}
\put(3,0){\line(1,0){3}} \multiput(0,0)(4,0){2}{\ten}
\put(-0.5,0.4){$-\beta_1$} \put(3.5,0.4){$-\beta_n$}
\put(6,0){\ten} \put(5.8,-0.4){$v_2$}
\put(5.5,-0.5){\framebox(2,1)} \put(-2.5,0){$\Gamma_2:$}
 \end{picture}
\end{center}

\subsection{Lemma}\label{sok}
{\em Consider the following resolution graph $\t\Gamma_2$ which
 satisfies $(1)$--$(4)$ below:
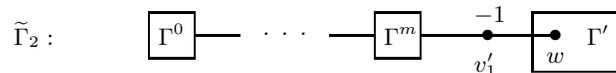
\begin{center}
\setlength{\unitlength}{0.6cm}
 \begin{picture}(11,2)(-1,-1)
\put(0,-0.5){\framebox(1,1)} \put(0.2,-0.2){$\Gamma^0$}
\put(1,0){\line(1,0){1}} \put(2.5,0){. . .}
\put(4,0){\line(1,0){1}} \put(5,-0.5){\framebox(1,1)}
\put(5.2,-0.2){$\Gamma^m$} \put(6,0){\line(1,0){3}}
\put(7.5,0){\ten} \put(7.2,0.3){$-1$} \put(7.2,-0.8){$v_1'$}
\put(8.5,-0.8){\framebox(2,1.3)} \put(9,0){\ten}
\put(8.8,-0.6){$w$} \put(9.7,-0.2){$\Gamma'$}
\put(-3,-0.2){$\t\Gamma_2:$}
 \end{picture}
\end{center}
\begin{enumerate}
 \item the subgraph on the left hand side of $w$ is a chain;
 \item $\Gamma^0$ consists of $(\alpha_0-1)$-vertices; the
       right-end has weight $-3$, all the others $-2$;
 \item $\Gamma^m$ consists of $(\alpha_m-2)$-vertices, all
       of them with weight $-2$;
 \item for $0<i<m$, $\Gamma^i$ consists of $(\alpha_i-2)$-vertices, the
       right-end has weight $-3$ and all the others $-2$.
\end{enumerate}

 Then $\t\Gamma_2$ is a negative definite plumbing
graph whose associated minimal graph is $\Gamma_2$. }

\begin{proof}
Glue the chain from the left hand side of $w$ with $\t\Gamma_0$,
blow down the $(-1)$--vertices and use the 0--absorption which
identifies \setlength{\unitlength}{0.6cm}
 \begin{picture}(5,0.8)(6,0)
\put(6.5,0){\line(1,0){2.5}} \put(7.5,0){\ten} \put(6.5,0){\ten}
\put(7.2,0.3){$-\beta_0$} \put(6.5,0.3){$0$}
\put(8.5,-0.5){\framebox(2,1)} \put(9,0){\ten}
\put(9.5,-0.2){$\Gamma_2$}
 \end{picture}
with $\Gamma_2$.
\end{proof}

\subsection{Corollary}\label{l:D2}
{\em  There exists a (non-minimal) resolution  graph, namely
$\t\Gamma_2$, which represents the splice diagram $\Delta_2$, it
has $\Gamma'$ as a subgraph, and supports a splice type
singularity.}

\begin{proof} We only have to show that
$\t\Gamma_2$ satisfies the monomial condition. 
Let $u$ be a node of $\Gamma_2$, let 
$C''$ a branch of it which contains the distinguished chain  involved in the 
splicing (otherwise the condition is trivial). Let $C$ be the branch of $u$
in $\Gamma$ which contains $C''$. Then there exists an effective cycle 
$A$ supported on $C$ so that $E_u^*+A$ is monomial in $\Gamma$.
Then the restriction $A|\Gamma'$ to $\Gamma'$ has the property that 
 $E_u^*(\Gamma')+A|\Gamma'$ is monomial and has negative intersection with 
$E_w$. Then $A|\Gamma'$ can be extended easily on the chain 
considered in (\ref{sok}) to an effective cycle $B$ on $\t\Gamma_2$
such that $E_u^*(\t\Gamma_2)+B$ is monomial in $\t\Gamma_2$.
\end{proof}

Let $H_{\t\Gamma_2,v_1'}(t)$ be the Hilbert series of associated
graded ring of the filtration with respect to $E_{v_1'}$ of
$\t\Gamma_2$. Then, again  by (\ref{t:pg}),  we have:

\subsection{Proposition}\label{p:2}
 $p_g(\Gamma_2)=H_{\t\Gamma_2,v_1'}|_{pc}+p_g(\Gamma')$.

\subsection{}
Finally we  consider the star--shaped graph $\Gamma_1$
(which automatically satisfies the monomial condition), hence 
(again by \thmref{t:pg}):

\subsection{Proposition}\label{p:3}
 $p_g(\Gamma_1)=H_{\Gamma_1,v_1}|_{pc}$.

\subsection{}\label{additivity}
 As a  consequence, from the above three propositions one gets:
\begin{equation}\label{addi}
 p_g(\Gamma)-p_g(\Gamma_1)-p_g(\Gamma_2)
=(H_{\Gamma,v_1}-H_{\t\Gamma_2,v_1'}-H_{\Gamma_1,v_1})|_{pc}. 
\end{equation}
The main point in this `additivity' formula is the following. 
Two terms in the expression  $\sigma+8p_g+c(\Gamma)=0$
(cf. \ref{4.1}) satisfy some additivity properties with respect to the
splicing. Namely, the additivity $\lambda(\Gamma)=\lambda(\Gamma_1)+
\lambda(\Gamma_2)$ (here we identify a plumbing graph with its 
plumbed 3-manifold) was proved independently by Akbulut-McCarthy, Boyer-Nicas
and Fukuhara-Maruyama (according to \cite{BoNi}). 
Let $b_1(G_i)$ be the first Betti number of the 
fiber $G_i$ of the fibered knot determined by $\Gamma_i$
and the splicing data ($i=1,2$). Then, by \cite[(6.4)]{NWuj}
(cf. also with \cite[(5.20)]{[51]}), one has 
$c(\Gamma)=c(\Gamma_1)+ c(\Gamma_2)-2b_1(G_1)b_2(G_2)$.
In particular, in order to run the induction of our proof, 
 we only need to verify 
$$p_g(\Gamma)=p_g(\Gamma_1)+p_g(\Gamma_2)+b_1(G_1)b_1(G_2)/4,$$
a fact already noticed in \cite[(6.3)]{NWuj}.
Therefore the following proposition, together with equation (\ref{addi}),
implies the theorem of the introduction.

\subsection{Proposition}\label{p:add2}
$(H_{\Gamma,v_1}-H_{\t\Gamma_2,v_1'}-H_{\Gamma_1,v_1})|_{pc}
=b_1(G_1)b_1(G_2)/4.$

\begin{proof} Note that in (\ref{p:add2}),  each periodic constant 
independently is a very complicated (Dedekind sum) expression. Nevertheless,
their combination provides the simple expression from the right hand side.
The proof is based on an explicit computation of the  Hilbert series
using (\ref{p:hp}).  The corresponding entries $m_{vw}$ might be determined 
from the splice diagrams 
using identities of \cite[\S 10]{EN}, or \cite[(9.1)]{NWuj}).
Note that the `unreduced' splice diagram associated with
$\t\Gamma_2$, having $v_1'$ as a node, is the following:

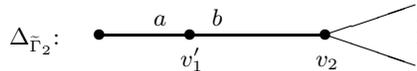
\begin{center}
\setlength{\unitlength}{0.6cm}
 \begin{picture}(11,2)(-6,-1)
\put(-2,0){\line(1,0){3}} \put(1,0){\line(3,1){2}}
\put(1,0){\line(3,-1){2}} \put(3,-0.2){$\vdots$} \put(-2,0){\ten}
\put(-2.2,-0.7){$v_1'$} \put(-1.5,0.2){$b$} \put(1,0){\ten}
\put(0.8,-0.7){$v_2$} \put(-2,0){\line(-2,0){2}} \put(-4,0){\ten}
\put(-2.8,0.2){$a$} \put(-6,-0.2){$\Delta_{\t\Gamma_2}$:}
 \end{picture}
\end{center}

We will use the notation of (\ref{FIG}) and let $\cV_2$ 
be the set of vertices of $\Gamma_2$. 
Then  $m_{v_1v_1}=ab$ and $m_{v_1w_i}=ab/a_i$.
Consider the integers  $p_i:=m_{v_1w_i}/b=a/a_i$ and $q_w:=m_{v_1w}/a$
for all $1 \le i \le r$ and  $w \in \cV_2$ respectively. 
Finally, set $g(t)=\prod_{w \in \cV_2}(1-t^{q_{w}})^{\delta_w-2}$. 
Then 
\begin{equation}\label{eq:H1}
 H_{\Gamma,v_1}(t)=\frac{g(t^a)(1-t^{ab})^{r-1}}{\prod_{i=1}^r(1-t^{p_ib})}.
\end{equation}
\no
Since for any $w \in \cV_2$,  $m_{v_1w}$ in ${\Gamma}$
and $m_{v_1'w}$ in  ${\t\Gamma_2}$ are the same, one has:
\begin{equation}\label{eq:H21'}
 H_{\t\Gamma_2,v_1'}(t)=\frac{g(t^a)}{1-t^b}.
\end{equation}
Clearly
\begin{equation}\label{eq:H11}
 H_{\Gamma_1,v_1}(t)
=\frac{(1-t^{ab})^{r-1}}{(1-t^a)\prod_{i=1}^r(1-t^{p_ib})}.
\end{equation}
We define the functions $Q_1(t)$ and $Q_2(t)$ by
$$
Q_1(t)=\frac{(1-t^{a})^{r-1}}{\prod_{i=1}^r(1-t^{p_i})}
-\frac{1}{1-t}, \quad Q_2(t)=g(t)-\frac{1}{1-t}.
$$
Then from \eqref{eq:H1}, \eqref{eq:H21'} and \eqref{eq:H11} we
have
\begin{equation}\label{eq:HQ}
   H_{\Gamma,v_1}(t)-
  H_{\t\Gamma_2,v_1'}(t)-H_{\Gamma_1,v_1}(t) =Q_1(t^b)Q_2(t^a)
-\frac{1}{(1-t^a)(1-t^b)}.
\end{equation}

On the other hand, consider the characteristic polynomial 
$P_i(t)=\det(I-th_{1,i})$ (where $I$ denotes the identity matrix)
of the monodromy $h_{1,i}\:H^1(G_i) \to H^1(G_i)$ 
of the fibered knot with Milnor fiber $G_i$ ($i=1,2$). 
Using A'Campo's formula \cite{AC2}, or \cite[\S 11--12]{EN} one has
\begin{align*}
P_1(t)&=(1-t)\cdot
 \frac{(1-t^a)^{r-1}}{\prod_{i=1}^r(1-t^{p_i})}, \\
P_2(t)&=(1-t)g(t).
\end{align*}
It is known for fibered links in integral homology 3--spheres
(or using A'Campo's formula is easy to check) 
that $P_i(1)=1$ and the derivative also satisfies $P'_i(1)=b_1(G_i)/2$.

Notice also that $Q_i(t)=(P_i(t)-1)/(1-t)$, hence $Q_i(1)=-P'_i(1)$. 
Therefore, since 
$Q_i$'s are polynomials, $(Q_1(t^b)Q_2(t^a))|_{pc}=Q_1(1)Q_2(1)=
b_1(G_1)b_1(G_2)/4$.
We end the proof with the remark that 
 $(1-t^a)^{-1}(1-t^b)^{-1}|_{pc}=0$ (since $gcd(a,b)=1$). 
\end{proof}

\end{document}